\documentclass{article}

\usepackage{graphicx}
\usepackage{amsmath}
\usepackage{amssymb}
\usepackage{jneurosci}
\usepackage[utf8]{inputenc}

\newcommand{\Ito}{It{\^o}}
\newcommand{\ZZ}{{\mathbf{Z}}}
\newcommand{\RR}{{\mathbf{R}}}
\newcommand{\PP}{{\mathbf{P}}}
\newcommand{\QQ}{{\mathbf{Q}}}
\newcommand{\cF} {{\mbox{$\cal F$}}}
\newcommand{\qed}{\rule{2mm}{2mm}}
\newcommand{\Lto}{\mbox{${\cal L}_2$}}
\newcommand{\transpose}{^\top}
\newcommand{\trace}{\mbox{\rm tr}}
\newcommand{\delvis}[2]{\frac{\partial #1}{\partial #2}}
\newcommand{\indicator}{\mathbf{1}}

\newcounter{envcounter} 
\newcommand{\myenv}[1]{\par\refstepcounter{envcounter}{\bf #1
    \theenvcounter:} }
\newenvironment{bemark}{\myenv{Remark}}{\mbox{ }\hfill$\qed$\par}
\newenvironment{example}{\myenv{Example}}{\mbox{ }\hfill$\qed$\par}

\begin{document}

\title{Transition probabilities for stochastic differential equations using the Laplace approximation: Analysis of the continuous-time limit}
\date{\today}
\author{Uffe Høgsbro Thygesen}
\date{DTU Compute, \texttt{uhth@dtu.dk}\\
  Technical University of Denmark\\
  DK-2800 Kongens Lyngby, Denmark}
\maketitle

\begin{abstract}
   We recently proposed a method for estimation of states and parameters in stochastic differential equations, which included intermediate time points between observations and used the Laplace approximation to integrate out these intermediate states. In this paper, we establish a Laplace approximation for the transition probabilities in the continuous-time limit where the computational time step between intermediate states vanishes. Our technique views the driving Brownian motion as a control, casts the problem as one of minimum effort control between two states, and employs a Girsanov shift of probability measure as well as a weak noise approximation to obtain the Laplace approximation. We demonstrate the technique with examples; one where the approximation is exact due to a property of coordinate transforms, and one where contributions from non-near paths impair the approximation. We assess the order of discrete-time scheme, and demonstrate the Strang splitting leads to higher order and accuracy than Euler-type discretization. Finally, we investigate numerically how the accuracy of the approximation depends on the noise intensity and the length of the time interval.
\end{abstract}

\section{Introduction}
\label{sec:intro}

Computation of transition probabilities is a fundamental problem for numerical analysis of stochastic differential equations, and arises in many different contexts. Here, our primary motivation comes from time series analysis, i.e., estimating states and parameters based on incomplete and/or noisy measurements of states in discrete time. Then, the joint density of the states are determined by the transition probabilities of the underlying stochastic differential equation. Recently \cite{Thygesen2025sdeA}, we considered an approach to this problem based on inserting unobserved states at intermediate time points between observations, and using the Laplace approximation to integrate out unobserved random variables; the present paper is a continuation of this research.

For general stochastic differential equations, the transition probabilities cannot be found in closed form, but a range of different numerical methods can be applied for their approximation \cite{Kloeden1999}. The probabilities are governed by the forward Kolmogorov equation (a.k.a. the Fokker-Planck equation), and numerical analysis of this partial differential equation is feasible with standard techniques (e.g., finite volume) in low-dimensional state spaces; up to three, say. Such a discretization of space can lead to approximating the entire time series analysis problem with a Hidden Markov Model \cite{Thygesen2023sde}, which is particularly useful for models with very strong nonlinear effects such as multimodal posteriors \cite{Thygesen2009b}. However, grid-based methods are infeasible for high-dimensional state spaces and inefficient when the posterior probability is centered on a small region of state space.  

When times between measurements are small, approximate expressions for the transition probabilities can be derived from time discretization schemes \cite{Kloeden1999}. Most basically, the explicit Euler-Maruyama method for an \Ito\ equation leads to a Gaussian approximation for the transition probability. This approximation can serve as a basis for time series analysis which was noted already in the 1970's \cite[e.g.]{PrakasaRao1979} and elaborated on many times since; for example, the approach can be refined with improved estimates for the conditional mean and variance \cite{Kessler1997}. Increased accuracy can be obtained with stronger discretization schemes \cite{Pilipovic2024} which makes the approach applicable to situations with larger time steps.

Approximations to the transition probabilities can also be obtained with linearization techiques \cite{Jazwinski1970}. Here, a nominal trajectory is determined as the solution to an ordinary differential equation, while deviations from this nominal trajectory are approximated with a Gaussian process which is determined from the linearization of the dynamics around the nominal path.  This is particular useful when nonlinearities in the drift function are important, also at the time scale between measurements, while the stochastic noise is weak. In the context of time series analysis, this leads to the Extended Kalman Filter, which have many variants \cite{Jimenez2003,Simon2006}. 

In this paper, we focus on yet a different approach to obtaining approximating transition densities, viz., the Laplace approximation. This approximation assists analysis of models with unobserved random variables, which are to be integrated out, by approximating the integrand with a Gaussian for which the integral is available in closed form. In the context of time series analysis, this approach has the appeal that it can also address the uncertainty of the states at the times of measurements, and thereby provides a unifying framework. \citetext{Karimi2014} investigated the application of the Laplace approximation to problems of filtering and estimation, under the assumption of additive noise; see also \cite{Karimi2016} and \cite{Thygesen2023sde}. The main contribution of \cite{Thygesen2025sdeA} was to extend this methodology beyond additive noise. 

The contribution of the present paper is to examine the continuous-time limit of the discrete-time schemes considered in \cite{Thygesen2025sdeA} for approximating transition densities. That is, we investigate the limit where the computational time steps between intermediate time points vanishes. An analysis with a similar objective, but somewhat different methods, has previously been performed by \citetext{Markussen2009}. Here, our analysis (section \ref{sec:continuous-laplace}) consists of four main elements: First, we identify the so-called \emph{most probable path} connecting the start and end points, by means of a minimum effort optimal control problem. Next, we shift the underlying probability measure using a Girsanov transformation to center the probability at this most probable path. Then, we analyze the properties of the control problem in the neighborhood of the most probable path, using linearization techniques and, in particular, a certain Riccati equation from optimal control. Finally, we examine the variance of weak-noise fluctuations around the most probable path, using a Lyapunov equation. These ingredients can then be combined to yield a Laplace approximation to the path integral, i.e., the integral over sample space. 

At the end, our method allows the computation of approximated transition probabilities, by solving ordinary differential equations. The primary application of this method, we believe, is not so much computation of actual transition probabilities in specific time series analysis, but rather that the continuous-time limit facilitates the analysis of discrete-time algorithms, and in particular, to disentangle errors from the Laplace approximation and from the time discretization. To demonstrate this potential, we examine the accuracy of the Laplace approximation (section \ref{sec:examples}) and compare the continuous-time limit with discrete-time approximations derived both from an Euler-type algorithm and from one involving Strang splitting. We also examine the error in the Laplace approximation itself. 

\section{Transition densities and their discretized Laplace approximation}
\label{sec:discr-time-lapl}

We consider the Stratonovich stochastic differential equation  in $\RR^n$
\begin{equation}
  \label{eq:sde}
  dX_t = f(X_t) ~dt + g(X_t) \circ dB_t , \quad X_0=x_0, 
\end{equation}
where $\{ B_t : t \geq 0\}$ is $n$-dimensional Brownian motion on a probability space $( \Omega, \cF,\PP)$ and  $\circ dB_t$ indicates the Stratonovich integral. We assume throughout that $f$ and $g$ satisfy the standard conditions (see, e.g., \citenoparens{Thygesen2023sde}) for existence and uniqueness of a strong solution $\{X_t : t \geq 0 \}$, and that $g(x)$ is invertible for all $x$. This implies that $\{X_t\}$ will be a Markov process for which the transition densities 
\[
  p(0,x_0,T,x_T)
\]
are well defined. Our main objective is to approximate these densities numerically.

We now summarize our approach \cite{Thygesen2025sdeA} for establishing a Laplace approximation for the transition probability $p(0,x_0,T,x_T)$. We first considered the implicit centered Euler-type approximation 
\begin{equation}
  \label{eq:Strat-approximation}
  X_{t+h} - X_t = \frac 12 (f(X_{t}) + f(X_{t+h})) ~h  + \frac 12 (g(X_{t}) + g(X_{t+h})) ~(B_{t+h} - B_t) . 
\end{equation}
Note that the noise intensity $g$ is evaluated both at $X_t$ and $X_{t+h}$, for consistency with the Stratonovich interpretation. Based on this discretization, \citetext{Thygesen2025sdeA} established the following approximation to the transition probabilities, valid for small time steps $h$:
\begin{multline}
  p(0,x,h,y) \approx \\ \left| I - \frac h2 \nabla f(y) - \sum_k \frac {b_k}2  \nabla g_k(y)  \right| \cdot \left| \frac{g(x) + g(y)}2\right|^{-1}  | 2 \pi h|^{-n/2} e^{- |b|^2/(2h)} .
  \label{eq:p-Euler}
\end{multline}
Here, $b=(b_1,\ldots,b_n)$  is the Brownian increment $B_h$ which brings the system from state $X_0=x$ to $X_h=y$ using the approximation \eqref{eq:Strat-approximation}, i.e., 
\begin{equation}
  \label{eq:b-formula}
  b = (g(x) + g(y))^{-1} (2 y - 2 x  - (f(x) + f(y)) h) . 
\end{equation}
To approximate the transition density $p(0,x_0,T,x_T)$ over large time intervals $[0,T]$, \citetext{Thygesen2025sdeA} next divided the entire time interval into equidistant subintervals using a partition $\{0,h,2h,\ldots, Nh\}$ with $Nh=T$, with the corresponding discretized trajectory
\[
  \bar x = \{ x_0,x_h,x_{2h}, \ldots, x_T\} . 
\]
Here, the end points $x_0$ and $x_T$ are fixed, while the remaining points in the discretized track are found by optimization as described in the following. To this discretized trajectory corresponds a sequence of increments of Brownian motion, say $\{ b_1,b_2,\ldots, b_N\}$, where each increment $b_i\in \RR^n$ is found from \eqref{eq:b-formula}, with $x=x_{(i-1)h}$, $y=x_{ ih}$. 
With this, we can find the joint log-density of the increments $\{ b_i\}$, as a function of the discretized track: 
\[
  \psi(\bar x )=  \sum_{i=1}^N [ \log \phi( b_i / \sqrt{h})  - \frac n2 \log h ]
  = - \frac{Nn}2 \log (2 \pi h) - \sum_{i=1}^N  \frac 12 \frac{|b_i|^2}h   . 
\]
Here, $\phi(\cdot)$ is the p.d.f. of a standard multivariate Gaussian random variable. We now find the most probable path
\[
  \bar x^* = \{ x_0, x_h^*,x_{2h}^* ,\ldots, x_T \} = \mbox{Arg}\max_{\bar x} \psi(\bar x) \quad \mbox{s.t. } x_0, ~x_T \mbox{ fixed,}
\]
assuming that there is indeed a unique global optimum. Combining the expression for the short-term transition density \eqref{eq:p-Euler} with the Laplace approximation, integrating over the Brownian increments, \citetext{Thygesen2025sdeA} then established the Laplace approximation
\begin{multline}
  \label{eq:Laplace-Stratonovich}
  \hat p(0,x_0,T,x_T ; h) = \\
  | H / (2 \pi) |^{-1/2} \exp(-\psi(\bar x^*)) \cdot \prod_{i=1}^{N} \frac{\left| I - \frac h2 \nabla f(x^*_{i})  - \sum_k \frac {b^*_{k(i-1)}}2  \nabla g_k (x^*_{i})\right| }
  {\left|\frac 12 (g(x^*_i) + g(x^*_{i-1}))\right|}
\end{multline}
where $H$ is the Hessian of $\psi$ at $\bar x^*$, and $b^*_{ki}$ is the corresponding increment of the $k$'th component of the Brownian motion over the $i$'th time step. This approximation arises when using the Laplace approximation to the integral over all increments of Brownian motion, i.e., over $\{ b_{ki} : ~k=1,\ldots,n, ~i = 1,\ldots, N\}$, subject to the constraint that these increments bring the system from $x_0$ to $x_T$. \citetext{Thygesen2025sdeA} demonstrated that this approximation was straightforward to implement and gave good approximations to the transition densities $p(\cdot)$ with modest computational effort. This, in turn allows estimation of states and parameters in multivariate stochastic differential equations with non-additive noise based on noisy and partial discrete-time observations of the state.

\section{The continuous-time limit}
\label{sec:continuous-laplace}

The Laplace approximation $\hat p(0,x_0,T,x_T; h)$ of $p(0,x_0,T,x_T)$ in  \eqref{eq:Laplace-Stratonovich} is based on a time discretization $h$, so the end result is affected by errors from both this discretization as well as from the Laplace approximation itself. It is useful to consider the limiting case where the time step $h$ vanishes, to separate these error sources, and thus to elucidate the error in the Laplace approximation. The continuous-time limit may also be useful for deriving stronger time discretizations. For these reasons, we examine here the continuous-time limit $h\rightarrow 0$ of the approximation $\hat p(0,x_0,T,x_T; h)$ from \eqref{eq:Laplace-Stratonovich}, keeping $x_0$, $x_T$ and $T$ fixed.

The outline of the argument, which establishes the continuous-time limit, is as follows: We first establish the ``most probable path'' (section \ref{sec:mpt}), around which the Laplace approximation is formed. We do this by considering the optimal control problem of how the noise can take the system from $x_0$ to $x_T$ with the least possible effort. Then, we shift the underlying probability measure (section \ref{sec:girsanovs-theorem}), to center the probability at this most probable path. This requires the state-feedback solution to the optimal control problem, viz. Bellman's value function, and  involves Girsanov's theorem, which establishes the Radon-Nikodyn derivative between the two measures. We can then establish the Laplace approximation by means of a weak-noise approximation (section \ref{sec:laplace-vs.-weak}); this requires a linearization of the optimal control problem, and next an analysis of how weak noise builds variance around the most probable path (section \ref{sec:small-noise-expans}). Thus, the computational vehicles involve a differential Riccati equation for the optimal control problem and a differential Lyapunov equation for the variance. Combining these elements, we reach the final continuous-time Laplace approximation \eqref{eq:final-result} in section \ref{sec:small-noise-expans}.

\subsection{The most probable path}
\label{sec:mpt}

We first examine how the mode $\bar x^*$ of the state sequence behaves as the time discretization $h$ becomes fine. Maximizing $\psi(\bar x)$ is equivalent to minimizing the criterion 
\begin{equation}
  \label{eq:crit}
  \frac12\sum_{i=0}^{N-1} \frac{| b_i|^2}{h} .
  \end{equation}
If we write the increments in the Brownian motion as $b_i = U_i h$, where $U_i \in \RR^n$,
we see that the optimization problem converges to the continuous-time problem of minimizing
\begin{equation}
  \label{eq:least-effort-problem}
  \min \int_0^T \frac 12 |U_t|^2 ~dt 
\end{equation}
where the minimization is over all trajectories $\{X_t: 0 \leq t \leq T\}$ and all driving noise terms $\{ U_t : t\in [0,T]\}$ which satisfy  the constraints
\begin{equation}
  \label{eq:state-constraints}
  dX_t = f(X_t) ~dt + g(X_t) U_t ~dt , \quad
  X_0=x_0, \quad X_T= x_T. 
\end{equation}
We call the solution $\{ (X_t,U_t) : 0 \leq t \leq T \}$ to this optimization problem the \emph{most probable path.} From the point of view of the stochastic differential equation (\ref{eq:sde}), this corresponds to finding the most probable realization of the Brownian motion which is consistent with the SDE and the two end points.
From the point of view of control theory, this is a \emph{minimum effort} control problem of steering a system from point $x_0$ to point $x_T$ in time $T$, where effort is understood in the sense of the \Lto-norm of the control $U:[0,T]\mapsto \RR^n$.
A substantial body of theory has been developed to analyze such control problems; see for example, the textbooks by  \citetext{BrysonHo1975,Bertsekas2005,Liberzon2012}. To recap this theory, we can pose the Lagrangian (a.k.a. the running cost)
\begin{equation}
  \label{eq:Lagrangian}
  L(x,u) = \frac 12 |u|^2
\end{equation}
so that the minimization problem~\eqref{eq:least-effort-problem} appears similar to the Least Action problems of analytical mechanics \cite{Lanczos70,Goldstein2002}. From this we form the pre-Hamiltonian $K$ (by several authors termed the Hamiltonian; note also differences in sign conventions) 
\[
  K(x,u,\lambda) = \lambda\transpose [ f(x) + g(x) u]  + L(x,u)
\]
where $\lambda$ is the co-state or the  ``shadow price'' which acts as a Lagrange multiplier on the constraint of the state dynamics in \eqref{eq:state-constraints}. In optimality, the pre-Hamiltonian is minimized w.r.t. $u$ at each point in time, so we define the optimal control law
\[
  \mu(x,\lambda) = \mbox{Arg} \min_u K(x,u,\lambda) ,
\]
where we assume that there is indeed a unique minimizing argument $u$ for each $(x,\lambda)$; this is the case for input-affine dynamics \eqref{eq:state-constraints} and convex quadratic Lagrangians \eqref{eq:Lagrangian}. This leads to the Hamiltonian (a.k.a. the optimized Hamiltonian or the control Hamiltonian)
\[
  H(x,\lambda) = \min_u K(x,u,\lambda)  . 
\]
In our case, with the Lagrangian~\eqref{eq:Lagrangian}, we get explicitly
\begin{equation}
  \label{eq:Hamiltonian}
  H(x,\lambda) = \lambda\transpose f(x) - \frac 12 | g\transpose(x) \lambda |^2 , \quad \mu(x,\lambda) = - g\transpose(x) \lambda. 
\end{equation}
With this Hamiltonian formulation, the optimal control law can be found by calculus of variations, and the resulting Euler-Lagrange equations take a form given by the celebrated Maximum Principle of Pontryagin \cite{Liberzon2012}. That is, let $\{(\bar X_t,\bar U_t) \}$ be the most probable path and control; then there exists a corresponding co-state $\{\bar \Lambda _t\}$ such that 
\begin{equation}
  \label{eq:Uopt}
  \bar U_s = \mu(\bar X_s,\bar \Lambda _s) = - g\transpose(\bar X_s) \bar \Lambda_s
\end{equation}
and such that the canonical equations 
\begin{equation}
  \label{eq:canonical}
 \dot {\bar X}_s = \delvis H\lambda ({\bar X}_s,\bar \Lambda_s) , \quad - \dot {\bar \Lambda}_s = \delvis Hx ({\bar X}_s,{\bar \Lambda}_s)  
\end{equation}
hold, as well as the boundary conditions $\bar X_0=x_0$, $\bar X_T=x_T$. These are necessary first-order conditions; we assume that there is a unique solution to these equations so that they characterize the most probable path uniquely, but return to this topic in the discussion. Instead of following Pontryagin, we can alternatively follow Bellman \cite{Bertsekas2005} and consider the cost-to-go $V(t,x)$ from an intermediate state $X_t=x$, defined by
\[
  V(t,x) = \inf \left\{ \frac 12 \int_t^T |U_s|^2 ~ds + \chi(X_T) : \dot X_s = f(X_s) + g(X_s) U_s, ~X_t = x \right\}
\]
where the infimum is over all paths and controls $\{ (X_s,U_s) : t \leq s \leq T\}$ which satisfy the given constraints. In terms of dynamic programming,  $V$ is the value function which satisfies the Hamilton-Jacobi-Bellman equation 
\begin{equation}
  \label{eq:HJB}
  \dot V + H(x,\nabla V) = 0 , \quad
  V(x,T) = \chi(x),
\end{equation}
at least provided $V$ is smooth, which we assume. The terminal cost $\chi(X_T)$ is yet unspecified but included so that we can manipulate the optimal path to end at $x_T$. To see how, note \cite{Liberzon2012} that the most probable path $\{\bar X_t\}$ and the associated co-state $\{\bar \Lambda_t\}$ can be seen as characteristics of this partial differential equation equation. Thus, we can identify the co-state $\bar \Lambda_t$ with the gradient $\nabla V(\bar X_t,t)$; in particular, the terminal cost $\chi$ must satisfy
\[
  \bar \Lambda_T = \nabla \chi(x_T)
\]
for the most probable path to end at $\bar X_T=x_T$. Apart from this constraint, $\chi(\cdot)$ is yet unspecified.

The identity $\bar \Lambda_t = \nabla V(\bar X_t,t)$ allows us to write the optimal control strategy \eqref{eq:Uopt} in state feedback form
\[
  \bar U_s = \mu( \bar X_s, \nabla V(\bar X_s,s))  = - g\transpose (\bar X_s) \nabla V(\bar X_s,s) . 
\]
Moreover, the framework allows us to evaluate the performance of not just the optimal path, but also every other path. Indeed, let $\{X_t,U_t\}$ be any path, which is not necessarily optimal, but which does satisfy the differential equation 
\[
  \dot X_t = f(X_t) + g(X_t) U_t ;
\]
with no constraints on initial and terminal point. Then a standard computation, which we  include for completeness,
shows that
\begin{align*}
  \frac 12 \int_0^T | U_t|^2 ~dt & = \int_0^T L(X_t,U_t) ~dt \\
                                                                         &= \int_0^T [ K - \nabla V \cdot \dot X ] ~dt \\
                                                                         & = \int_0^T [ K - H - \frac d{dt} V ] ~dt \\
  & = V(X_0,0) - V(X_T,T) + \int_0^T [ K-H] ~dt .
\end{align*}
Here, we have omitted the arguments of $K(X_t,U_t,\nabla V(X_t,t))$, $H(X_t,\nabla V(X_t,t))$, and $V(X_t,t)$. This result decomposes the cost into an optimal one which depends on the end points, and a deviation from optimality which can be seen as the integral of decision errors $K-H$ accumulated along the path. With the specific Lagrangian \eqref{eq:Lagrangian}, these decision errors can be written in terms of the control signal
\[
  K-H = \frac 12 | U_t - \mu(X_t,\nabla V(X_t,t)) |^2 , 
\]
and, with the dynamics \eqref{eq:state-constraints}, more explicitly in terms of the state trajectory,
\[
  K-H = \frac 12 \left| (g^{-1}) \left( \dot X_t - f + gg\transpose  \nabla V \right) \right|^2.
\]
Here we have omitted the arguments of $f(X_t)$, $g(X_t)$, and $\nabla V(X_t,t)$, for compactness. This characterization is useful for analysis of near-optimal paths.

\subsection{Centering the probability at the optimal path}
\label{sec:girsanovs-theorem}

We now compare the original system
\[
  dX_t = f(X_t) ~dt + g(X_t) \circ dB_t
\]
where $\{B_t\}$ is $\PP$-Brownian motion, with a controlled system
\begin{equation}
  \label{eq:controlled-system}
  dX_t = f(X_t) ~dt - g(X_t) g\transpose(X_t) \nabla V(X_t,t) ~dt + g(X_t) \circ dW_t .
\end{equation}
That is, we define
\[
  W_t = B_t + \int_0^t g\transpose (X_t) \nabla V(X_s,s) ~ds . 
\]
The point of this shift is that with $W_t\equiv0$, we obtain the most probable path connecting $X_0=x_0$ and $X_T=x_T$, as identified in the previous. 
We now introduce a new probability measure $\QQ$ on the sample space such that $\{W_t\}$ is $\QQ$-Brownian motion. By Girsanov's theorem \cite{Oeksendal2010}, we have
\begin{equation}
  \label{eq:Girsanov}
  \PP [Z] = \QQ \left[ Z  L_T \right]
\end{equation}
where $Z$ is a random variable and where $\PP[\cdot]$ and $\QQ[\cdot]$ denote expectation of the random variable in question w.r.t. the probability measures $\PP$ and $\QQ$, respectively. Here, $L_T$ is the likelihood, also referred to as the Radon-Nikodyn derivative of $\PP$ w.r.t. $\QQ$, which is in itself a random variable given by \cite{Oeksendal2010} 
\begin{align*}
  L_T
  & =\exp \int_0^T (  \frac 12 | \nabla V \cdot g  |^2 ~dt  + \nabla V  \cdot g ~ dB_t ) . 
\end{align*}
Note that the stochastic integral is an \Ito\ integral, even if our starting point was a Stratonovich equation.
To examine the exponent in $L_T$, we first use the Hamilton-Jacobi-Bellman equation~(\ref{eq:HJB}) with the Hamiltonian \eqref{eq:Hamiltonian} to rewrite the likelihood as 
\[
  L_T  =\exp \int_0^T (  \dot V ~dt + \nabla V \cdot f  ~dt  + \nabla V \cdot g ~dB_t ) 
\]
and next introduce the scalar process $\{Y_t : 0 \leq t \leq T \}$ given by $Y_t = V(X_t,t)$. By \Ito's formula,
\[
  dY_t = [ \dot V + \nabla V \cdot f_I + \frac 12 \trace (g\transpose \nabla \nabla V g )] ~dt + \nabla V \cdot g ~dB_t 
\]
where $f_I$ is the \Ito\ drift, given by $f_I(x) = f(x) + \frac 12 \sum_k \nabla g_k  \cdot g_k $. Therefore, we have
\[
d \log L_t = dY_t - \frac 12 \nabla V \cdot \sum_k \nabla g_k  \cdot g_k ~dt - \frac 12 \trace (g\transpose \nabla \nabla V g ) ~dt. 
\]
Defining the differential operators $M_k$ for $k=1,\ldots,n$ by $M_k V = \nabla V \cdot g_k$, we can use the identity
\[
  \sum_k M_k M_k V = \trace ( g\transpose \nabla \nabla V g) + \nabla V \cdot \sum_k \nabla g_k \cdot g_k
\]
to rewrite this relationship between $d \log L_t$ and $dY_t$ as
\[
  d\log L_T = dY_t - \frac 12 \sum_k M_k M_k V ~dt  
\]
and thus
\begin{equation}
  \label{eq:LT-result}
  L_T = \exp\left( Y_T - Y_0 -
    \frac 12 \int_0^T \sum_k M_k M_k V ~dt   \right) . 
\end{equation}
It may be illuminating to note that if $g$ is the constant identity matrix, then $\sum_k M_k M_k V = \Delta V$; in general, $\sum_k M_kM_k V$ combines curvature of $V$ and noise intensity $g$ in a coordinate-free way to yield a noise sensitivity. Thus, the likelihood $L_T$ decreases with the least required effort $V(X_0,0) - V(X_T,T)$ to steer the deterministic system from $X_0$ to $X_T$, and with this noise sensitivity $\sum_k M_kM_kV$ integrated along the realized path. 

Now, let $dx$ be a small set in state space which contains the end point $x_T$, and consider the random variable $Z=\indicator_{X_T \in dx }$, the indicator variable that the realized path ends near $x_T$. Then, combining \eqref{eq:Girsanov} and \eqref{eq:LT-result}, we obtain
\begin{multline}
  \label{eq:exact-transform}
  \PP [ \indicator_{X_T \in dx } ] =
  \\
  \exp(V(x_T,T)-V(x_0,0))  \QQ \left[ \indicator_{X_T \in dx}
  \exp\left( - \frac 12 \int_0^T \sum_k M_k M_k V  ~dt \right) 
  \right] . 
\end{multline}

Here, we have used that $dx$ is a small neighborhood where $V(\cdot,T)$ can be considered constant and equal to $V(x_T,T)$. We emphasize that up to this point, no other approximations have been made. The aim of the following sections is to apply a Laplace approximation to the $\QQ$-expectation in \eqref{eq:exact-transform}. 

\subsection{Laplace vs. weak noise approximations}
\label{sec:laplace-vs.-weak}

We now note a basic connection between Laplace approximations and weak noise approximations of probability densities. This connection seems natural, but we have not been able to find it precisely stated in the literature. We first state this result in the finite-dimensional context: Consider an $n$-dimensional Gaussian random variable $Z\sim N(0,\Sigma)$ with  $\Sigma >0$, and define 
\[
  X = h(Z)
\]
where $h:\RR^n\mapsto \RR^m$ is $C^2$. We take $h(0)=0$ for simplicity, and assume that $n\geq m$ and that $H:=\nabla h(0)$ is surjective, viewing $\nabla h$ as an $m$-by-$n$ matrix. Our aim is to approximate the probability density of $X$ evaluated at $x=0$, comparing two approximation techniques. We first consider an approximation based on linearization and weak noise, as in the well-known delta method. Specifically, define
\[
  X_\epsilon := \epsilon^{-1} h(\epsilon Z)
\]
then as $\epsilon \rightarrow 0$, $X_\epsilon$ converges in distribution to a Gaussian with mean 0 and variance $H \Sigma H\transpose$. Therefore, the p.d.f. of $X_\epsilon$ at 0 converges to
\begin{equation}
  \label{eq:weak-noise-density}
  | 2 \pi H\Sigma H\transpose  |^{-1/2} . 
\end{equation}
Note that this expression is exact if $m=n$ and $h$ is a diffeomorphism, but our interest is in the case $n>m$. On the other hand, we can write the density of $X$ at $x=0$ as the following limit:
\[
  f_X(0) = \lim_{\delta\rightarrow 0} \int_\ZZ f_Z(z) ~\phi( h(z) / \delta ) ~\delta^{-m} ~ dz  . 
\]
Let $\hat f_X(0; \delta)$ denote the Laplace approximation of this integral. To find $\hat f_X(0;\delta)$, we first evaluate the integrand at $z=0$ to get
\[
  | 2 \pi \Sigma |^{-1/2} |2 \pi | ^{-m/2} \delta^{-m}  , 
\]
and next find the Hessian of the negative log-integrand at $z=0$:
\[
  \Sigma^{-1} + \delta^{-2} H\transpose H . 
\]
Therefore, the Laplace approximation of the integral is 
\[
  \hat f_X(0; \delta) = | 2 \pi \Sigma |^{-1/2} |2 \pi | ^{-m/2} \delta^{-m} |   \Sigma^{-1} + \delta^{-2} H\transpose  H |^{-1/2} (2\pi)^{n/2} 
\]
which reduces to 
\[
 \hat f_X(0;\delta ) =  |2 \pi | ^{-m/2} \delta^{-m} |   I + \delta^{-2}  \Sigma^{1/2} H\transpose H \Sigma^{1/2} |^{-1/2} . 
\]
The determinant can be written
\[
  |   I + \delta^{-2}  \Sigma^{1/2} H\transpose H\Sigma^{1/2} |
  = \prod_{i=1}^m (1+\delta^{-2} \sigma_i^2)
\]
where $\{ \sigma_1,\ldots,\sigma_m \}$ are the non-zero singular values of $H\Sigma^{1/2}$.
In the limit $\delta \rightarrow0$, we get
\[
  |   I + \delta^{-2}  \Sigma^{1/2} H\transpose H\Sigma^{1/2} | \cdot \delta^{2m} \rightarrow | H \Sigma H\transpose| 
\]
and therefore
\begin{equation}
  \label{eq:Laplace-density}
  \hat f_X(0 ; \delta ) \rightarrow |   2 \pi H\Sigma H\transpose    |^{-1/2} .
\end{equation}
Comparing \eqref{eq:weak-noise-density} and \eqref{eq:Laplace-density}, we see the Laplace approximation of the p.d.f. of $X$ at $x=0$ can be obtained also by the weak noise approximation.

The previous was for a finite-dimensional random variable $Z$, but the result carries over to the case where $X$ is the solution to a stochastic differential equation at a given point of time, say $X_T$. To see this, note that $X_T$ can be approximated using a time discretization \cite{Kloeden1999} in which case $Z$ will contain the discrete-time increments of the driving Brownian motion. Note that the simplicity of the result requires that we evaluate the density at the point $x$ which is obtained with $Z=0$; otherwise, second-order terms will appear, as shown by \citetext{Markussen2009} in a somewhat different formulation. This simplicity is the motivation for the change of measure in section \ref{sec:girsanovs-theorem}: The Laplace approximation of the transition density under $\QQ$, say, $q(0,x_0,T,x_T)$, can be obtained by weak noise approximation, i.e., local linearization in the sense of \citetext{Jazwinski1970} and \citetext{Jimenez2003}, since $x_T$ is the endpoint of the noise-free ($W\equiv 0$) trajectory. 
  
\subsection{Small noise expansion}
\label{sec:small-noise-expans}

The result of the previous section leads us to consider a weak noise approximation of the $\QQ$-expectation in   \eqref{eq:exact-transform}, where we replace the $\QQ$-Brownian motion $\{W_t\}$ with $\{\epsilon W_t\}$. This requires a local analysis near the optimal path, involving the quadratic approximation of the value function $V(x,t)$ near the optimal path:
\[
  V(x,t) = V({\bar X}_t,t ) + \bar \Lambda_t\transpose (x-{\bar X}_t) + \frac 12 (x-{\bar X}_t)\transpose Q_t (x-{\bar X}_t) + o(|x-{\bar X}_t|^2) . 
\]
Here, we have already established the constant and linear terms: 
\[
  V({\bar X}_t,t ) = \int_t^T \frac 12 | \bar U_s|^2 ~ds + \psi({\bar X}_T) 
\mbox{ 
  and
  }
  \nabla V({\bar X}_t,t) = \bar \Lambda_t ,
\]
and for the quadratic term, a standard result from optimal control \cite{BrysonHo1975} is that the Hessian $Q_t = \nabla \nabla V(\bar X_t,t)$ is governed by the Riccati equation 
\begin{equation}
  \label{eq:Riccati}
  \dot Q_t + H_{xx} + H_{x\lambda} Q_t + Q_t H_{\lambda x} + Q_t H_{\lambda \lambda} Q_t = 0, \quad Q_T = \nabla \nabla \psi({\bar X}_T) .
\end{equation}
Here, the derivatives of the Hamiltonian $H(x,\lambda)$ are evaluated along the optimal path $(\bar X_t, \bar \Lambda_t)$. For the scalar case $n=1$, we find 
\begin{align*}
  H_{xx}  &= \lambda f'' - \frac 12 \lambda ^2 (g^2)'' , \\
  H_{x \lambda } &= f' - \lambda (g^2)', \\
  H_{\lambda \lambda } & = - g^2  .
\end{align*}

\begin{example}
  For systems with affine drift and additive noise, $f(x) = Ax + b$, $g(x)=G$, we have $H_{xx}=0$, so a solution is $Q_t \equiv 0$. One interpretation of this observation is that the curvature of the value function is needed for nonlinear systems, because the effect of the noise depends on the trajectory, whereas for linear systems superposition applies.
\end{example}

\begin{example}
  For general systems $(f,g)$, if the end point $x_T$ is that obtained with $B_t \equiv 0$, then $\bar \Lambda_t \equiv 0$, and $H_{xx}=0$. So also in this case, a solution is $Q_t \equiv 0$, reflecting that the entire step of shifting the probability measure is unnecessary.
\end{example}

We are now ready to make a weak noise approximation of the  $\QQ$-expectation in \eqref{eq:exact-transform}, where we replace $W_t$ with $\epsilon W_t$. To this end, note first that when $\epsilon \rightarrow  0 $, we have $X_t \rightarrow \bar X_t$ w.p. 1 (in particular, $X_T \rightarrow x_T$), and 
\begin{equation}
  \label{eq:iMMVdt-limit}
  \int_0^T M_k M_k V(X_t,t) ~dt \rightarrow \int_0^T ( g_ k \cdot Q_t g_k + \bar \Lambda_t\transpose \cdot \nabla g_k ~g_k) ~dt .
\end{equation}
Here, the integral on the right hand side is along the optimal path, i.e., $g_k$ and its Jacobian is evalauted at $\bar X_t$. Next, to first order in $\epsilon$, $X_t$ is Gaussian distributed with expectation $\bar X_t$ and a variance $\epsilon^2 \Sigma_t$ given by the linearization of system dynamics around the most probable path $\{ {\bar X}_t\}$, i.e., the Lyapunov equation
\begin{equation}
  \label{eq:lyap}
  \dot \Sigma_t = A_t \Sigma_t + \Sigma_t A_t\transpose + G_tG\transpose_t , \quad \Sigma_0=0.  
\end{equation}
Here, $G_t = g(\bar X_t)$ while $A_t$ is the Jacobian of the drift field $f - gg\transpose \nabla V$ evaluated at $(\bar X_t,t)$; in the scalar case, this is $A_t = f' - g^2 Q_t - 2gg' \bar \Lambda_t$. Combining, we get
\begin{align*}
  & \QQ \left[ \indicator_{X_T \in dx}
  \exp\left(\frac 12 \int_0^T  \sum_k M_kM_k V   ~dt \right) 
    \right] \\
  & \approx |2 \pi\epsilon^2 \Sigma_T|^{-1/2} \exp \left( -\frac 12 \int_0^T ( \trace g\transpose Q_t g + \bar \Lambda_t \cdot \sum_k \nabla g_k g_k  \right) ~|dx|.
\end{align*}
According the principle of weak noise approximation, we insert $\epsilon=1$ in this expression. We next combine this with \eqref{eq:exact-transform} to obtain an approximate expression for $\PP [ \indicator_{X_T \in dx } ]$, which then leads to 
\begin{multline}
  \hat p(0,x_0,T,x_T) = 
  \\
    |2 \pi \Sigma_T|^{-1/2} \exp \left( -\frac 12 \int_0^T ( |{\bar U}_t|^2 + \trace g\transpose Q_t g + \bar \Lambda_t \cdot \sum_k \nabla g_k g_k ) ~dt \right). 
  \label{eq:final-result}
\end{multline}
This is the final result for the Laplace approximation of the transition density $p(0,x_0,T,x_T)$. To summarize, its evaluation requires that we first solve the canonical equations \eqref{eq:canonical} with the boundary conditions $\bar X_0=x_0$, $\bar X_T = x_T$. Then, we compute the optimal control \eqref{eq:Uopt}, and solve the Riccati equation \eqref{eq:Riccati} as well as the Lyapunov equation \eqref{eq:lyap}. Finally, we evaluate the integrals in \eqref{eq:final-result}. 

\begin{bemark}
  [Coordinate transformations]
  \label{remark:coordinates}
  The resulting Laplace approximation for the transition probabilities of $X_T$ is independent of the choice of coordinate systems. That is, if we define the transformed process $\{ Z_t=\eta (X_t) \}$ where $\eta(\cdot)$ is a diffeomorphism, so that the transition density $q$ of $\{Z_t\}$ satisfies the usual coordinate transformation $q(0,z_0,T,z_T) \cdot | \nabla \eta(x_T)| = p (0,x_0,T,x_T) $, then the same relationship holds for the Laplace approximations. To see this, we note first that according to the differential geometric approach to optimal control, the Most Probable Path is independent of the coordinate system; $\bar Z_t = \eta(\bar X_t)$, and this path is obtained with the same minimum effort control signal $\{ \bar U _t\}$ in the two coordinate systems. The value function also transforms, i.e., if the value function in $x$-coordinates is $V(x,t)$, then the value function is $z$-coordinates is $V(\eta^{-1}(z),t)$; the same transform applies to $M_kM_k V$. Since coordinate transformations acts identically in standard deterministic calculus and in Stratonovich calculus, also $g\transpose Q_t g$ and $\bar \Lambda_t \nabla g_k g_k$ are invariant under coordinate transformations, so that this also applies to the entire exponent in \eqref{eq:final-result}. Finally, the covariance $\Sigma_t$ of the linearized system transforms to
  $\nabla \eta(\bar X_t) \Sigma_t \nabla \eta\transpose(\bar X_t)$. Thus, when we change coordinates, the only change to the expression \eqref{eq:final-result} is the inclusion of the term $|\nabla \eta (\bar X_T)|$.
\end{bemark}

\section{Examples}
\label{sec:examples}

\subsection{Geometric Brownian motion}
\label{sec:geom-brown-moti}

For Geometric Brownian Motion, we show in this example that the Laplace approximation is, in fact, exact. This is to be expected in view of remark \ref{remark:coordinates}, since the process can be transformed (with a log transform) into Brownian motion with drift, for which the Laplace approximation is exact. Therefore, the example serves mainly to verify and showcase the machinery. The starting point is the Stratonovich equation

\begin{equation}
  \label{eq:GBM-Stratonovich}
  dX_t = rX_t ~dt + \sigma X_t \circ dB_t ,
\end{equation}
for which we know that the transition probabilities are log-Gaussian; specifically, $\log X_t | X_0 = x \sim N(\log x + rt,\sigma^2 t)$. For the minimum-effort control problem \eqref{eq:least-effort-problem}, the pre-Hamiltonian is 
\[
  K(x,u,\lambda) = \lambda x (r+\sigma u) + \frac 12 u^2 
\]
  and stationarity  $\partial K/\partial u = 0$ gives the optimal strategy and the Hamiltonian:

\[
  \mu(x,\lambda) = -\sigma \lambda x , \quad
    H(x,\lambda) = \lambda x r - \frac 12 \sigma^2 \lambda^2 x^2 . 
\]

The Hamiltonian is constant along optimal paths \cite{BrysonHo1975},
so
\[
  \bar X_t \bar \Lambda_t = \alpha 
\]
must hold for some constant $\alpha$. Therefore, the optimal control is constant in time, $\bar U_t=-\sigma \alpha$, and the optimal state trajectory is exponential in time and given by the linear equation:
\[
  \dot {\bar X}_t = \bar X_t ( r- \alpha \sigma^2 )
    . 
\]
The  parameter $\alpha$ must be chosen so that $\bar X_T=x_T$, i.e., 
\[
  (r-\alpha \sigma^2) T = \log (x_T/x_0) 
\]
which leads to the most probable path 
\[
  \bar X_t = x_T^{t/T} x_0^{1-t/T} . 
\]
Thus, the most probable path interpolates the end points geometrically, as expected. The value of $\alpha$ also defines the co-state:
\[
  \bar \Lambda_t = \sigma^{-2} (r - T^{-1} \log (x_T/x_0)) x_T^{-t/T} x_0^{t/T-1} .
\]
In particular, at time $T$, we get
\[
  \bar \Lambda_T = \sigma^{-2}(r - T^{-1} \log (x_T/x_0)) / x_T.
\]
To identify the value function, we guess a solution of the form
\[
  V(x,t) = \beta t + \gamma \log x . 
\]
To match the boundary condition $V'(x_T,T) = \bar\Lambda_T$, we get
\[
  \gamma / x_T = \sigma^{-2} (r - T^{-1} \log (x_T/x_0)) / x_T \Leftrightarrow \gamma  = \sigma^{-2} (r - T^{-1} \log (x_T/x_0))  . 
\]
Note that only $V'(x_T,T)$ is specified; the particular choice of a logarithm in $V(\cdot,T)$ is for simplicity.
In turn, inserting the candidate solution $V$  into the HJB equation, we get
\[
  \beta +  \gamma r - \frac 12 \sigma^2 \gamma^2 = 0
\]
This fixes $\beta = \frac 12 \sigma^2 \gamma^2 - \gamma r$, and with these value of $\beta$ and $\gamma$, our candidate $V$ does indeed satisfy the HJB equation. 
With this value function, we find 
\[
  MV(x) = V' \sigma x = \gamma \sigma
\]
and therefore $MMV (x) = 0$. Also this is as expected. Therefore, the expression \eqref{eq:exact-transform} simplifies to
\[
  \PP [ \indicator_{X_T\in dy} ] = \exp( \beta T + \gamma (\log y - \log x))   \QQ [ \indicator_{X_T\in dy} ] .
\]
With this $V$,  the optimally controlled system \eqref{eq:controlled-system} becomes
\[
  dX_t = rX_t ~dt - \sigma^2 \gamma X_t ~dt  + \sigma X_t \circ dW_t . 
\]
Thus, the optimally controlled process is geometric Brownian motion, just as the original process; the difference being that with $W_t \equiv 0$ we obtain $X_T=x_T$. As an intermediate checkpoint, we can therefore verify the correctness of \eqref{eq:exact-transform}, using the known transition probabilities of geometric Brownian motion. We have, from the p.d.f. of the log-Gaussian distribution
\[
  \PP [ \indicator_{X_T\in dx} ] =
  x_T^{-1} (2 \pi \sigma^2 T)^{-1/2} \exp(\frac 12 \frac{(\log x_T - \log x_0 - rT)^2}{\sigma^2 T} ~dx ,
\]
and similarly, 
\[
  \QQ [ \indicator_{X_T\in dx} ] =
  x_T^{-1} (2 \pi \sigma^2 T)^{-1/2} .
\]
To verify \eqref{eq:exact-transform}, we must therefore show that
\[
  \exp\left( \beta T + \gamma \log \frac {x_T}{x_0} \right) = \exp\left(-\frac 12 \frac{(\log x_T - \log x_0 - rT)^2}{\sigma^2 T}\right). 
\]
Inserting $\beta = \frac 12 \sigma^2 \gamma^2 - \gamma r$, and using $\gamma  = \sigma^{-2} (r - T^{-1} \log (y/x))$, we get
\[
  \frac 12 \sigma^2 \gamma^2 T - \gamma r T + \gamma \log \frac yx = - \frac 12  \gamma^2 \sigma^2 T
\]
which again reduces to $0=0$, thus verifying the claim. This, of course, just confirms that we have transformed the measures correctly, and we are only able to conduct this verification because we know the transition probabilities both under $\PP$ and $\QQ$.

To conduct the small noise approximation, we first note from the function $V(x,t)$ that
\[
  Q_t = V''(\bar X_t,t) = - \gamma \bar X_t^{-2}
\]
and therefore
\[
  \dot Q_t = 2 \gamma \bar X_t^{-3} \dot {\bar X}_t .
\]
We now check that this Hessian satisfies the Riccati equation \eqref{eq:Riccati}, the terms in which become
\[
  H_{xx} = -\gamma^2 \bar X_t^{-2} \sigma^2 , \quad H_{x\lambda} = r -  2\gamma \sigma^2 , \quad
  H_{\lambda \lambda} = - \sigma^2 \bar X_t^2
\]
(here, we have evaluated the derivatives of the Hamiltonian along the path), so that the Riccati equation itself becomes (using $\bar X_t^{-1} \dot {\bar X}_t = T^{-1} \log (x_T/x_0)$)
\[
  2 \gamma \bar X_t^{-2} T^{-1} \log \frac {x_T}{x_0}    -\gamma^2 \bar X_t^{-2} \sigma^2 - 2(r-2\gamma \sigma^2) \gamma \bar X_t^{-2} - \gamma^2 \sigma^2 \bar X_t^{-2} = 0
\]
which simplifies to
\[
  2 T^{-1} \log \frac {x_T}{x_0} - \gamma \sigma^2 - 2(r-2\gamma \sigma^2) - \gamma \sigma^2 = 0
\]
or
\[
  T^{-1} \log \frac {x_T}{x_0} + \gamma \sigma^2  - r  = 0
\]
which is true. We have therefore confirmed that the Hessian $Q_t$ of the value function $V(\cdot,t)$ at $\bar X_t$ satisfies the Riccati equation.

To find the transition probabilities, we linearize the system dynamics around the most probable path: We have
  \[
    A_t = f'(\bar X_t) + g'(\bar X_t) \bar U_t = r + \sigma \bar U_t = T^{-1} \log (x_T/x_0), \quad
    G_t = \sigma {\bar X}_t = \sigma x_T^{t/T} x_0^{1-t/T} . 
  \]
  It is convenient to rewrite the Lyapunov equation \eqref{eq:lyap} in terms of a relative variance $R_t$ given by 
  \[
    S_t = R_t ({\bar X}_t)^2
  \]
  
  so that the relative variance $R_t$ satisfies 
  \[
    \dot R_t = \frac{\dot S_t}{({\bar X}_t)^2} - 2\frac {S_t \dot {\bar X}_t}{({\bar X}_t)^3} = 2 T^{-1} \log\frac {x_T}{x_0} R_t + \sigma^2 - 2 R_t T^{-1} \log \frac {x_T}{x_0} = \sigma^2 . 
  \]
  Here we have used $\dot {\bar X} = T^{-1} \log (x_T/x_0) ~{\bar X}_t$. Thus
  \[
    R_t = \sigma^2 t
  \]
  and, in particular,
  \[
    S_T = \sigma^2 x_T^2 T . 
  \]
  Thus, the Laplace approximation of the p.d.f. of $X_T$ at $x_T$ is
  \[
    \hat p(0,x0,T,x_T) = \frac1{\sqrt{2 \pi S_T}} e^{-W(x_T,T)}
    =
    \frac1{\sqrt{2 \pi \sigma^2 x_T^2 T}} \exp\left( - \frac 12 
      \frac{ \log \frac {x_T}{x_0} - rT}{\sigma^2 T}
    \right)
  \]
  which we recognize as the p.d.f. of a log-Gaussian random variable $X_T$ with log-mean $\log x_0 + rT$ and log-variance $\sigma^2 T$.  This is consistent with 
  \[
    X_T = x_0 e^{rT + \sigma B_T} 
  \]
  which is the solution of the Stratonovich equation~\eqref{eq:GBM-Stratonovich}. We see that in this case, the Laplace approximation is in fact exact, $\hat p(\cdot) = p(\cdot)$. To reiterate, we knew this on forehand, given that the the problem can be log-transformed to Brownian motion with drift, for which the transition probabilities are Gaussian and therefore the Laplace approximation is exact, and we have included the calculations to demonstrate the details of the machinery. In more typical situations, where we do not know the answer on forehand, this analysis will largely by replaced by numerical solutions of the relevant ordinary differential equations.

\subsection{The double well}
\label{sec:double-well}

This example shows an inevitable shortcoming of the Laplace approximation, using the well-known stochastic double well system
\[
  dX_t = [ X_t - X_t^3] ~dt + \sigma \circ dB_t . 
\]
We consider the transition probability from $x_0=0$ to $x_T=0$. According to the Laplace approximation, we have
  $Q_t\equiv0$ as the optimal control is $\bar U_t \equiv 0$. The variance of the linearization is given by the Lyapunov equation
  \[
    \dot \Sigma_t = 2 \Sigma_t + 1
  \]
  so $\Sigma_t = \frac 12 (e^{ 2 t}-1)$ and  
  \[
    \hat p(0,0,t,0) = | 2 \pi \Sigma_t|^{-1/2} . 
  \]
  In particular, the approximation $\hat p(0,0,t,0)$ vanishes as $t\rightarrow\infty$, in contrast to the true transition probability, which converges to that given by the stationary distribution. The reason for this discrepancy is that the true transition probability involves also contributions from paths which are far from the most probable path; for example, paths that transit to one of the wells, fluctuates there, and then returns to the origin at time $T$.

\subsection{Order analysis for discretization schemes}
\label{sec:order-analysis-cox}

Here, we examine by means of a numerical experiment how quickly the discrete-time Laplace approximation given in section \ref{sec:discr-time-lapl} converges to the continuous-time limit established in section \ref{sec:continuous-laplace}, i.e., the errors between the discrete expression \eqref{eq:Laplace-Stratonovich} and the continuous expression  \eqref{eq:final-result}. To this end, we consider the Cox-Ingersoll-Ross process, which is governed by the \Ito-sense stochastic differental equation
\begin{equation}
  \label{eq:CIR}
  dX_t = \lambda(\xi - X_t) ~dt + \gamma \sqrt{X_t} ~dB_t . 
\end{equation}
We take parameters $\lambda=1$, $\xi=1$ and $\gamma=1/2$. We aim to compute the transition probability $p(0,x_0,T,x_T)$ with $x_0=0.75$,  $T=1$ and $x_T$ near 1.5. For the continuous-time limit of section \ref{sec:continuous-laplace}, we use the standard ODE solver \texttt{ode} from the \texttt{R} package \texttt{deSolve} with default settings to solve the Hamiltonian system \eqref{eq:canonical} for the most probable path. To address that the initial condition $\bar \Lambda_0$ is unknown, we use shooting, i.e., we guess on $\bar \Lambda_0$ and refine the guess until we have $\bar X_T \approx 1.5$ with sufficient accuracy. We find $\bar \Lambda_0 = -2.106$ which results in $\bar X_T \approx 1.500024$, we take this value as $x_T$, so that the following analysis is not affected by shooting error. We then solve the Riccati equation \eqref{eq:Riccati} as a terminal value problem using the terminal condition $Q_T=0$, and next the Lyapunov equation \eqref{eq:lyap}. The coefficients in these equations depend on the trajectory $\bar X_t$, $\bar \Lambda_t$; this is stored with a (rather excessive) time resolution of 0.001 and interpolated linearly. The details of the implementation can be inspected in the accompanying code repository \cite{Thygesen_SDETMB}. We find $\hat p(0,x_0,T,x_T) = 0.256$. For this Cox-Ingersoll-Ross process, the true transition probabilities are well known \cite[exercise 9.8]{Thygesen2023sde}; see also \cite{Thygesen_SDEtools}. We have $p(0,x_0,T,x_T)\approx 0.257$, so that the absolute error on the continuous-time Laplace approximation is about $10^{-3}$. 

\begin{figure}
  \centering
  \includegraphics[width=0.7\textwidth]{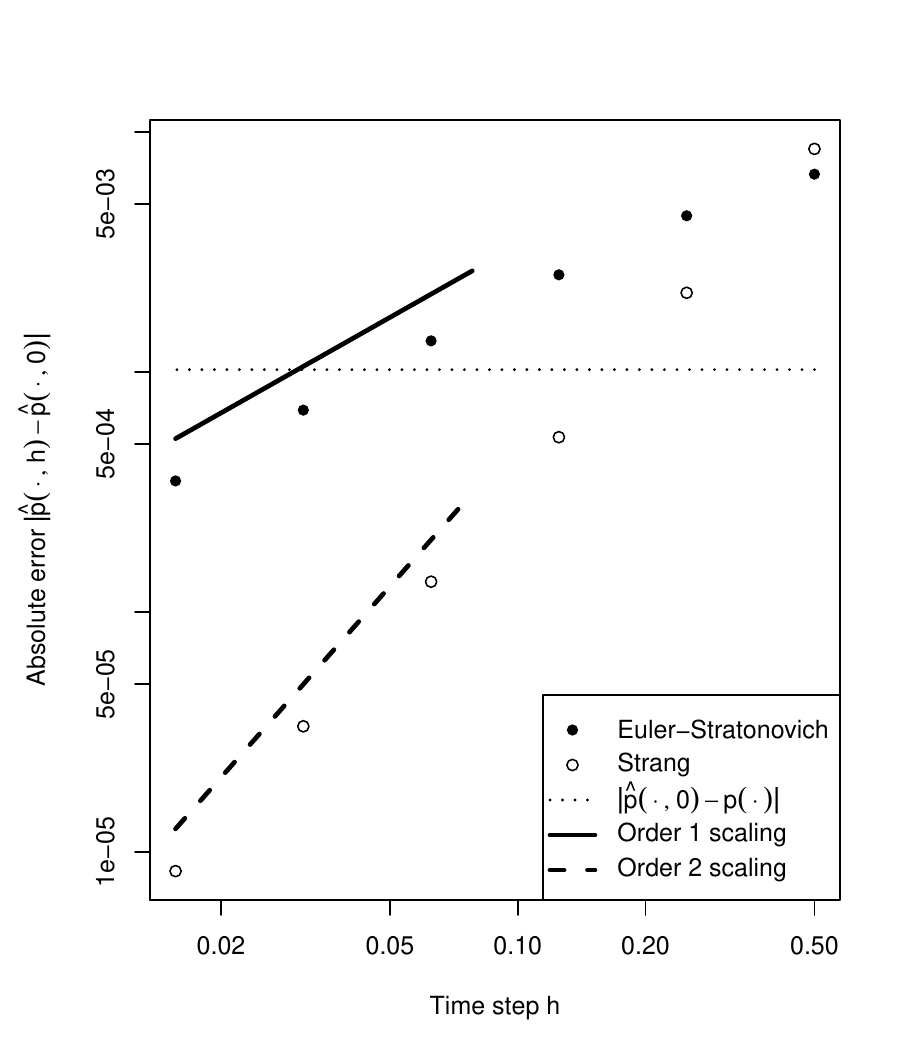}
  \caption{Errors in computing the transition densities for the Cox-Ingersoll-Ross process: Comparison of the discrete-time Laplace approximation $\hat p(\cdot,h)$ for different time steps $h$ and for the two discretization methods , the limiting continuous-time Laplace approximation $\hat p(\cdot,0)$, and the true transition density $p(\cdot)$. See text for full description.}
  \label{fig:CIR-transprob}
\end{figure}

Next, we compute the same transition probabilities using the discrete-time Laplace approximation of section \ref{sec:discr-time-lapl}, i.e., the expression \eqref{eq:Laplace-Stratonovich}, for different choices of the time step $h$. The results are seen in figure \ref{fig:CIR-transprob}, which shows the absolute difference $|\hat p(\cdot,h) - \hat p(\cdot,0)|$ between the discrete-time Laplace approximation transition probability using a time step $h$, and the continuous-time Laplace approximation, corresponding to a time step $h=0$. The graph also includes the error $|p(\cdot) - \hat p(\cdot,0)|$ between the true transition probability $p(\cdot)$ and the continuous-time limit $\hat p(\cdot,0)$. For reference, the plot also includes a straight line with slope 1, corresponding to a numerical order 1 (note that the plot is double logarithmic). 

The plot indicates that the discretization error is of order 1, and that - for these parameters - the time step should be around 0.05 or smaller, if one aims to have a discretization error which is no greater than the error from the Laplace approximation itself. 

\subsubsection{A Strang splitting algorithm}
\label{sec:strang-splitt-algor}

One of our motivations for pursuing the continuous-time limit of section \ref{sec:continuous-laplace} is that it facilitates assessment of different discretization methods; this is the first requirement for developing improved numerical methods. To demonstrate this, we derive and demonstrate in this section an improved discretization method for the Cox-Ingersoll-Process based on Strang splitting; see \cite{Foster2024,Pilipovic2024} and the references therein.

To recap, the numerical method of Strang splitting was originally developed for ordinary differential equations but has recently become popular for stochastic differential equations. It is based on a decomposition of the underlying stochastic differential equation into two terms. The decomposition can be done in different ways which gives rise to different methods; here, we write the stochastic differential equation \eqref{eq:CIR} as
\[
  dX_t = dX_t^{(1)} + dX_t^{(2)} \mbox{ where } dX^{(1)}_t = \lambda (\xi - X_t) ~dt , \quad dX^{(2)}_t = \gamma \sqrt{X_t} \circ dB_t .
\]
The advantage of this decomposition is that the two flows corresponding to the two terms have analytical solutions. Specifically,
\[
  dX_t = dX^{(1)}_t \Leftrightarrow X_t = \xi + (X_0 - \xi) e^{-\lambda t},
\]
while
\[
  dX_t = dX^{(2)}_t \Leftrightarrow X_t = (\sqrt{X_0} + \frac \gamma 2 B_t )^2 .
\]
In the following, we make use of these analytical solutions; in the discussion we will address the question what could have been done if these analytical solutions were not available. Strang splitting now updates the state by first following the one flow for half a time step $h/2$, then the other for a full time step $h$, and finally the first flow for half a time step  $h/2$ again. This principle is also found in leapfrogging and symplectic integration. Specifically, letting $X^{(1)}$ and $X^{(2)}$ be intermediate variables, we get for the Cox-Ingersoll-Ross process
\begin{align*}
  X^{(1)} &= \xi + (X_0 - \xi) e^{-\lambda h/2}, \\
  X^{(2)} &= (\sqrt{X^{(1)}} + \frac \gamma 2 B_h )^2 , \\
  X_h &= \xi + (X^{(2)} - \xi) e^{-\lambda h/2} . 
\end{align*}
Here, we use these approximate relationships to compute the conditional p.d.f. of $X_h$ at a given point $y$ conditional on $X_0=x$, i.e., the transition density $p(0,x,h,y)$. To this end, we follow the procedure of section \ref{sec:discr-time-lapl}; that is, given $X_0=x$ and $X_h=y$, we first solve for $X^{(1)}$, $X^{(2)}$, and then $B_h$. We then compute the p.d.f. of $B_h$ at the value we have found, and use the implicit function theorem to get the p.d.f. of  $X_h$ at $y$. After this, we use the same technique as in section \ref{sec:discr-time-lapl} to perform the discrete-time Laplace approximation; the final implementation only differs in a couple of lines of code, as is evident from the sample code accompanying this paper \cite{Thygesen_SDETMB}.

The results are seen in figure \ref{fig:CIR-transprob}. We see that Strang splitting in this way and for this example leads to second order accuracy. If the criterion is to have a discretization error which is smaller than that of the Laplace approximation, then Strang splitting requires roughly 4 times fewer time steps. While we have not done a detailed analysis of execution time, the total time to produce the data in this graph is around 0.1 seconds on a standard laptop, both for the Euler-Stratonovich method and for the Strang method. Thus, the improved accuracy of the Strang method for a given time step comes without any substantial cost in terms of execution time. Also, the coding complexity of the algorithms is roughly similar. In summary, we see no reason to prefer the Euler-Stratonovich algorithm over Strang splitting.

\subsection{The accuracy of the Laplace approximation}
\label{sec:accur-lapl-appr}

\begin{figure}
  \centering
  \includegraphics[width=\textwidth]{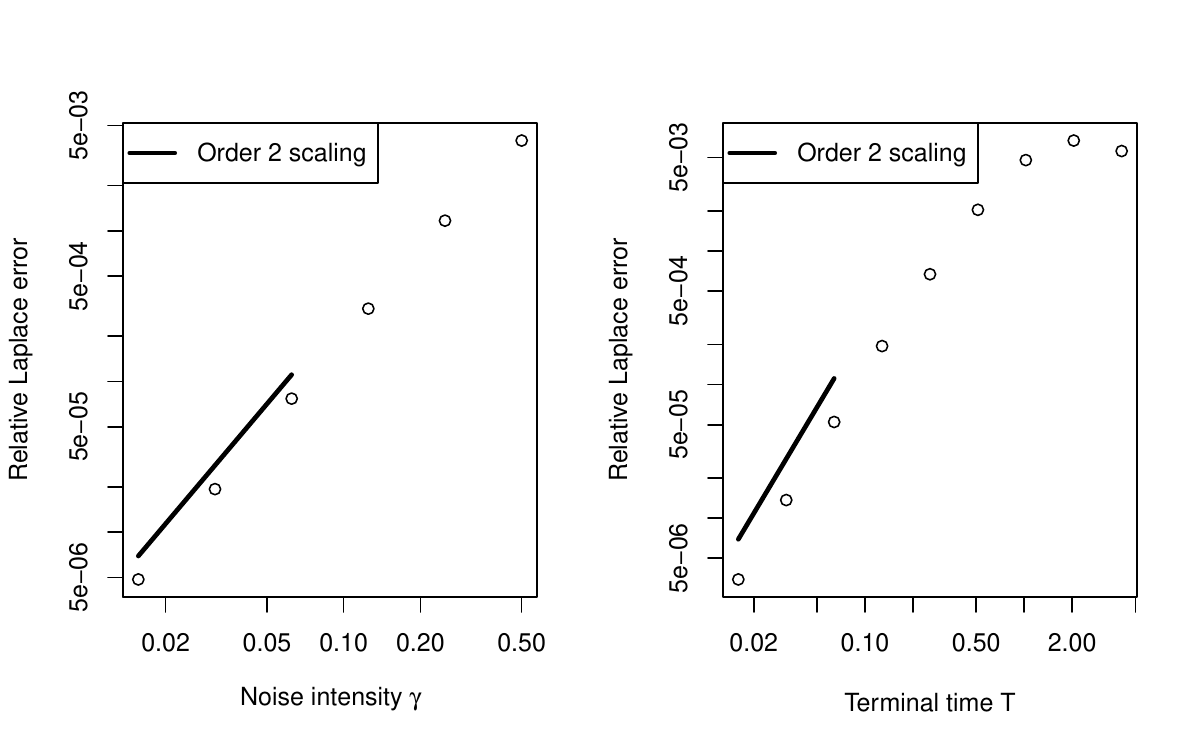}
  \caption{Errors in the continuous-time Laplace approximation for the Cox-Ingersoll-Ross process \eqref{eq:CIR}:
    Comparison of the theoretical transition densities with the expression \eqref{eq:final-result}.
    \emph{Left panel:} The effect of varying the noise intensity $\gamma$.
    \emph{Right panel: } The effect of varying the final time $T$.
    Base parameters are $\lambda=1$, $\xi=1$, $\gamma=0.5$, $T=1$, $x_0=0.75$.
    See the text for how the terminal point $x_T$ is chosen.}
  \label{fig:CIR-Laplace-error}
\end{figure}

We now investigate the accuracy of the continuous-time Laplace approximation \eqref{eq:final-result} for the Cox-Ingersoll-Ross process \eqref{eq:CIR}, and, in particular, how the relative error depends on the terminal time $T$ and the noise intensity $\gamma$. The results are seen in figure \ref{fig:CIR-Laplace-error} and indicate that the relative error scales with the square of $T$ for fixed $\gamma$, and with the square of $\gamma$ for fixed $T$. 

When conducting an analysis like this, a choice must be made how the end point $x_T$ is chosen. When varying the noise intensity $\gamma$, we kept the initial co-state $\Lambda_0$ constant, and solved the canonical equations to find the terminal point $x_T$. When varying the terminal time $T$, we varied the initial co-state as $\Lambda_0 = -2\exp(-T)$; this heuristic ensures that the terminal point $x_T$ lies approximately at the same quantile in the prior distribution of $X_T$.

\section{Discussion}
\label{sec:discussion-1}

Time series analysis based on stochastic differential equations require methods for computing transition densities. While we previously \cite{Thygesen2025sdeA} discussed discrete-time Laplace approximations to this end, we have in this paper pursued the continuous-time limit as the computational time step vanishes. This can also be seen as a Laplace approximation of the integral over sample space, i.e., path integrals. With our method, the computation of the Laplace approximation involves the solution of ordinary differential equations, viz., the canonical equations for the optimal path, the Riccati equation for the Hessian of the value function, and the Lyapunov equation for the variance of fluctuations around the most probable path. The method applies to quite general multivariate stochastic differential equations, including state-dependent noise intensity, and can easily be extended to also include explicit time dependence in the dynamics.

This paper has only done a first analysis of this problem, and a number of issues remain. For example, we have assumed that there is a unique optimal path, but it is easy to conceive singular situations where this assumption does not hold. From the theory of optimal control \cite{Liberzon2012}, we know that such so-called conjugate points can arise is some systems when the time interval becomes large enough, and can at least partly be diagnosed by finite escape times in the Riccati equation. A situation which is perhaps more critical for applications is that other locally optimal paths can contribute significantly to the total integral over path space, so that the Laplace approximation is inaccurate. We mentioned one such situation when discussing the double well in section \ref{sec:double-well}.

This leads to the broader question of accuracy of the Laplace approximation. We have demonstrated by example that the error increases with time $T$ and noise level $\gamma$ as $\gamma^2 T^2$. A theoretical underpinning of this results remains; this would also take into account the degree of non-linearity of the dynamics. For such studies, it is convenient to have benchmarks such as the Cox-Ingersoll-Ross process, where the transition densities are known while the Laplace approximation is not perfect; unfortunately, we are not aware of many such benchmarks. The motivation for further research in the accuracy of the Laplace approximation would include possible transformations and reformulations which could improve the accuracy, and it could also investigate the case where the noise intensity $g(x)$ is not square and invertible.

On a theoretical note, it would be interesting to investigate further the relationship between our results and those by \citetext{Markussen2009}. At first glance, our approaches are quite different, but both works are motivated by the same ambition and both arrive at ordinary differential equations for the Laplace approximation, focusing on the most probable element in the underlying Wiener space in a way which is independent of the coordinate system used in state space. Thus, it is plausible that the approaches are in fact parallel. Such a study could also include casting the limit of the discrete-time expression \eqref{eq:Laplace-Stratonovich} in terms of Fredholm determinants.  

The primary use of our apparatus is, we believe, conception and analysis of discretization techniques, whether based on numerical methods for stochastic differential equations as the Strang method in section \ref{sec:strang-splitt-algor}, or based on numerical solutions of the ordinary differential equations we have posed. Our numerical results indicate that the Euler-Stratonovich scheme leads to an order 1 while the Strang splitting leads to an order 2. These results are not surprising but we have not attempted a theoretical analysis to support them; note that standard results for the order of numerical schemes concern the initial value problem. As for adapting higher-order numerical methods for SDE's, note that our Laplace approximation grows in complexity if each computational time step and each dimension of state space involves more than one random variable.

Our application of Strang splitting in section \ref{sec:strang-splitt-algor} used that the two flows involved each had analytical solutions. The algorithm would have been very similar, though, if we had used a numerical solution of the pure-drift flow. An obvious candidate is Heun's method, as this gives the required second order accuracy, but also stiff solvers and higher-order methods could be useful, in particular for advection-dominated problems. This could extend the Strang splitting method to any model, where the pure noise part has a known analytical solution; this covers many models used in practice. 

In summary, our study adds theoretical support under the computational approach in \cite{Thygesen2025sdeA} by establishing a continuous-time limit. This allows a more specific analysis of approximation errors in discrete-time algorithms, by separating the Laplace error and the discretization error. In combination, this will facilitate time series analysis based on stochastic differential equations with non-additive noise, where the time between observations is not small compared to the time scales of the system dynamics. 

\bibliographystyle{jneurosci}
\bibliography{new}

\begin{thebibliography}{}

\bibitem[\protect\citeauthoryear{Bertsekas}{2005}]{Bertsekas2005}
Bertsekas DP (2005)
\newblock {\em Dynamic Programming and Optimal Control}, Vol.~1
\newblock Athena Scientific, Belmont, Massachusetts.

\bibitem[\protect\citeauthoryear{Bryson and Ho}{1975}]{BrysonHo1975}
Bryson AE, Ho YC (1975)
\newblock {\em Applied optimal control}
\newblock Taylor \& Francis.

\bibitem[\protect\citeauthoryear{Foster \bgroup et al.\egroup
  }{2024}]{Foster2024}
Foster JM, dos Reis G, Strange C (2024)
\newblock High order splitting methods for sdes satisfying a commutativity
  condition.
\newblock {\em SIAM Journal on Numerical Analysis}~62:\mbox{500--532}.

\bibitem[\protect\citeauthoryear{Goldstein \bgroup et al.\egroup
  }{2002}]{Goldstein2002}
Goldstein H, Poole C, Safko J (2002)
\newblock Classical mechanics.

\bibitem[\protect\citeauthoryear{Jazwinski}{1970}]{Jazwinski1970}
Jazwinski A (1970)
\newblock {\em Stochastic processes and filtering theory}
\newblock Academic Press, New York.

\bibitem[\protect\citeauthoryear{Jimenez and Ozaki}{2003}]{Jimenez2003}
Jimenez J, Ozaki T (2003)
\newblock Local linearization filters for non-linear continuous-discrete state
  space models with multiplicative noise.
\newblock {\em International Journal of Control}~.

\bibitem[\protect\citeauthoryear{Karimi and McAuley}{2014}]{Karimi2014}
Karimi H, McAuley KB (2014)
\newblock A maximum-likelihood method for estimating parameters, stochastic
  disturbance intensities and measurement noise variances in nonlinear dynamic
  models with process disturbances.
\newblock {\em Computers \& Chemical Engineering}~67:\mbox{178--198}.

\bibitem[\protect\citeauthoryear{Karimi and McAuley}{2016}]{Karimi2016}
Karimi H, McAuley KB (2016)
\newblock Bayesian estimation in stochastic differential equation models via
  {Laplace} approximation.
\newblock {\em IFAC-PapersOnLine}~49:\mbox{1109--1114}
\newblock 11th IFAC Symposium on Dynamics and Control of Process Systems
  Including Biosystems DYCOPS-CAB 2016.

\bibitem[\protect\citeauthoryear{Kessler}{1997}]{Kessler1997}
Kessler M (1997)
\newblock Estimation of an ergodic diffusion from discrete observations.
\newblock {\em Scandinavian Journal of Statistics}~24:\mbox{211--229}.

\bibitem[\protect\citeauthoryear{Kloeden and Platen}{1999}]{Kloeden1999}
Kloeden P, Platen E (1999)
\newblock {\em Numerical Solution of Stochastic Differential Equations}
\newblock Springer, New York, third edition.

\bibitem[\protect\citeauthoryear{Lanczos}{1970}]{Lanczos70}
Lanczos C (1970)
\newblock {\em The Variational Principles of Mechanics}
\newblock Dover Publications, Inc., New York, fourth edition
\newblock The first edition was published in 1949.

\bibitem[\protect\citeauthoryear{Liberzon}{2012}]{Liberzon2012}
Liberzon D (2012)
\newblock {\em Calculus of variations and optimal control theory: a concise
  introduction}
\newblock Princeton University Press.

\bibitem[\protect\citeauthoryear{Markussen}{2009}]{Markussen2009}
Markussen B (2009)
\newblock Laplace approximation of transition densities posed as {Brownian}
  expectations.
\newblock {\em Stochastic Processes and their
  Applications}~119:\mbox{208--231}.

\bibitem[\protect\citeauthoryear{{\O}ksendal}{2010}]{Oeksendal2010}
{\O}ksendal B (2010)
\newblock {\em Stochastic Differential Equations - An Introduction with
  Applications}
\newblock Springer-Verlag, sixth edition.

\bibitem[\protect\citeauthoryear{Pilipovic \bgroup et al.\egroup
  }{2024}]{Pilipovic2024}
Pilipovic P, Samson A, Ditlevsen S (2024)
\newblock Parameter estimation in nonlinear multivariate stochastic
  differential equations based on splitting schemes.
\newblock {\em The Annals of Statistics}~52:\mbox{842--867}.

\bibitem[\protect\citeauthoryear{{Prakasa Rao} and
  Rubin}{1979}]{PrakasaRao1979}
{Prakasa Rao} B, Rubin H (1979)
\newblock Asymptotic theory for process least squares estimators for diffusion
  processes
\newblock Technical report 79-13, Department of Statistics, Purdue University.

\bibitem[\protect\citeauthoryear{Simon}{2006}]{Simon2006}
Simon D (2006)
\newblock {\em Optimal State Estimation - {Kalman}, {$\cal H_\infty$}, and
  Nonlinear Approaches}
\newblock John Wiley \& Sons, Hoboken, New Jersey.

\bibitem[\protect\citeauthoryear{Thygesen \bgroup et al.\egroup
  }{2009}]{Thygesen2009b}
Thygesen UH, Pedersen MW, Madsen H (2009)
\newblock Geolocating fish using hidden {Markov} models and data storage tags
\newblock In Nielsen J, Arrizabalaga H, Fragoso N, Hobday A, Lutcavage M,
  Sibert J, editors, {\em Tagging and Tracking of Marine Animals with
  Electronic Devices}, Vol.~9 of {\em Reviews: Methods and Technologies in Fish
  Biology and Fisheries}, \mbox{pp. 277--293}. Springer.

\bibitem[\protect\citeauthoryear{Thygesen}{2023}]{Thygesen2023sde}
Thygesen UH (2023)
\newblock {\em Stochastic Differential Equations for Science and Engineering}
\newblock CRC Press/Taylor and Francis.

\bibitem[\protect\citeauthoryear{Thygesen}{2025a}]{Thygesen_SDETMB}
Thygesen UH (2025a)
\newblock {SDE-TMB}
\newblock R code available at github.com/Uffe-H-Thygesen/SDE-TMB.

\bibitem[\protect\citeauthoryear{Thygesen}{2025b}]{Thygesen_SDEtools}
Thygesen UH (2025b)
\newblock {SDEtools}
\newblock R package available at github.com/Uffe-H-Thygesen/SDEtools.

\bibitem[\protect\citeauthoryear{Thygesen and
  Kristensen}{2025}]{Thygesen2025sdeA}
Thygesen UH, Kristensen K (2025)
\newblock Inference in stochastic differential equations using the {Laplace}
  approximation: Demonstration and examples.
\newblock {\em arXiv}~\mbox{p. 2503.21358}
\newblock https://doi.org/10.48550/arXiv.2503.21358.

\end{thebibliography}

\end{document}